\documentclass[twoside,leqno]{article}  

\usepackage{ltexpprt} 
\usepackage{amssymb,amsmath}
\usepackage[ruled]{algorithm2e}
\usepackage{graphicx}
\usepackage{subfigure}
\usepackage{url}

\newcommand{\argmax}{\operatornamewithlimits{argmax}}
\newcommand{\argmin}{\operatornamewithlimits{argmin}}

\usepackage{etoolbox}
\newtoggle{use_ORTH_PNMF_instead_of_PNMF}
\toggletrue{use_ORTH_PNMF_instead_of_PNMF}
\iftoggle{use_ORTH_PNMF_instead_of_PNMF}
{
\newcommand{\ProjNMF}{O-PNMF}
\newcommand{\ProjNMFSVD}{O-P(SVD)}
}
{
\newcommand{\ProjNMF}{PNMF}
\newcommand{\ProjNMFSVD}{P(SVD)}
}

\begin{document}

\title{\Large Two Algorithms for Orthogonal Nonnegative Matrix Factorization 
\\ with Application to Clustering} 
\author{Filippo Pompili\thanks{Department of Electronic and Information Engineering, 
University of Perugia, Italy; \texttt{filippo.pompili@diei.unipg.it}.} 
\and 
Nicolas Gillis\thanks{Universit\'e de Mons, Department of Mathematics and Operational Research, Facult\'e Polytechnique, Mons, Belgium; 
\texttt{nicolas.gillis@umons.ac.be} (Corresponding author)} 
\and 
P.-A. Absil\thanks{Universit\'e catholique de Louvain, ICTEAM Institute, Avenue Georges Lemaitre 4, B-1348 Louvain-la-Neuve, Belgium; \texttt{pa.absil@uclouvain.be}}
\and 
Fran\c{c}ois Glineur${}^\ddagger$\thanks{Universit\'e catholique de Louvain, CORE, Voie du Roman Pays 34, B-1348 Louvain-la-Neuve, Belgium; \texttt{francois.glineur@uclouvain.be}.}}
\date{}

\maketitle

\begin{abstract} \small\baselineskip=9pt 
 Approximate matrix factorization techniques with both nonnegativity and orthogonality constraints, referred to as orthogonal nonnegative matrix factorization (ONMF), have been recently introduced and shown to work remarkably well for clustering tasks such as document classification. In this paper, we introduce two new methods to solve ONMF. First, we show mathematical equivalence between ONMF and a weighted variant of spherical $k$-means, from which we derive our first method, a simple EM-like algorithm. This also allows us to determine when ONMF should be preferred to $k$-means and spherical $k$-means. Our second method is based on an augmented Lagrangian approach. Standard ONMF algorithms typically enforce nonnegativity for their iterates while trying to achieve orthogonality at the limit (e.g., using a proper penalization term or a suitably chosen search direction). Our method works the opposite way: orthogonality is strictly imposed at each step while nonnegativity is asymptotically obtained, using a quadratic penalty. Finally, we show that the two proposed approaches compare favorably with standard ONMF algorithms on synthetic, text and image data sets.  
\end{abstract}


\section{Introduction}

We consider the orthogonal nonnegative matrix factorization (ONMF)
problem, which can be formulated as follows. Given an $m$-by-$n$
nonnegative matrix $M$ and a factorization rank $k$ (with $k <
n$), 
solve
\begin{subequations}  \label{eq:ONMF}
\begin{align}
 \min_{U \in \mathbb{R}^{m \times k}, V \in \mathbb{R}^{k \times n}} \;
& ||M-UV||_F^2  \label{eq:objfn}
\\   \text{subject to} \quad & U \geq 0, \ V \geq 0, \label{eq:nonneg}
\\  & VV^T = I_k,  \label{eq:orth}
\end{align}
\end{subequations}
where $\|\cdot\|_F$ denotes the Frobenius norm,~\eqref{eq:nonneg}
means that the entries of matrices $U$ and $V$ are nonnegative, and
$I_k$ stands for the $k\times k$ identity matrix. 

The ONMF problem~\eqref{eq:ONMF} can be viewed as the well-known nonnegative matrix
factorization (NMF) problem,~\eqref{eq:objfn}-\eqref{eq:nonneg}, with
an additional orthogonality constraint,~\eqref{eq:orth}, that
considerably modifies the nature of the problem. In particular, it is
readily seen that constraints~\eqref{eq:nonneg} and~\eqref{eq:orth}
imply that $V$ has at most one nonzero entry in each
column; we let $i_j$ denote the index of the nonzero entry (if any) in column $j$ of $V$.
Therefore, any solution $(U^*,V^*)$ of~\eqref{eq:ONMF} has the following property:
for $j=1,\ldots,n$, index $i_j$ is such that column $i_j$ of $U^*$ achieves the smallest angle with column $j$ of data matrix $M$, 
while $V^*(i_j,j)$ scales column $i_j$ of $U^*$ to make it as close as possible to column $j$ of $M$ (in the sense of the Euclidean norm). 
Hence it is clear that the ONMF problem relates to data clustering and, indeed, empirical evidence suggests that the
additional orthogonality constraint~\eqref{eq:orth} can improve
clustering performance compared to standard NMF or
$k$-means~\cite{Choi08algorithmsfor,Yang10linearand}.

Current approaches to ONMF problems are based on suitable
modifications of the algorithms developed for the original NMF
problem. They enforce nonnegativity of the iterates at each step, and
strive to attain orthogonality at the limit (but never attain exactly
orthogonal solutions). This can be done using a proper penalization
term~\cite{Ding06orthogonalnonnegative}, a projection matrix
formulation~\cite{Yang10linearand} or by choosing a suitable search
direction~\cite{Choi08algorithmsfor}. Note that, for a given data
matrix $M$, different methods may converge to different pairs
$(U,V)$, where the objective function~\eqref{eq:objfn} may take
different values.  
Furthermore, under random initialization, which is used by most NMF algorithms~\cite{Bou08SVD}, two runs of the same method may yield different results. This situation is due to the multimodal nature of the ONMF
problem~\eqref{eq:ONMF}---it may have multiple local minima---along
with the inability of practical methods to guarantee more than
convergence to local, possibly nonglobal, minimizers. Hence, ONMF methods not only differ in
their computational cost, but also in the quality of the clustering
encoded in the returned pair $(U,V)$ for a given problem.

In this paper, we first show the equivalence of ONMF with a weighted variant of spherical $k$-means, which leads us to design an EM-like algorithm for ONMF. We also explain in which situations ONMF should be preferred to $k$-means and spherical $k$-means. 
Then, we propose a new ONMF method, dubbed ONP-MF, that
relies on a strategy reversal: instead of enforcing nonnegativity of
the iterates at each step and striving to attain orthogonality at the
limit, ONP-MF enforces orthogonality of its iterates while obtaining
nonnegativity at the limit. 
A resulting advantage of ONP-MF is that rows of factor $V$ can be initialized directly with the right singular vectors of $M$ (which is the optimal solution of the problem without the nonnegativity constraints), whereas the other methods require a prior alteration of the singular vectors to make them nonnegative~\cite{Bou08SVD}. We show that, on
some clustering problems, the new algorithm outperforms other
clustering methods, including ONMF-based methods, in terms of
clustering quality.

The paper is organized as follows. In Section~\ref{equiv}, we analyze the relationship between ONMF and clustering problems and show that it is closely related to spherical $k$-means. Based on this analysis, we develop an EM-like algorithm which features a rank-one NMF problem at its core. This also allows us to shed some light on the differences between $k$-means, spherical $k$-means and ONMF, which we illustrate on synthetic data sets. 
Section~\ref{npomf} introduces another algorithm to perform ONMF using an augmented Lagrangian and a projected gradient scheme, which enforce orthogonality at each step while obtaining nonnegativity at the limit. Finally, in Section~\ref{nr}, we experimentally show that our two new approaches perform competitively with standard ONMF algorithms on text data sets and on different image decomposition problems. \\ 

This paper is an extended version of the proceedings paper \cite{PGAG}. 


\section{Equivalence of ONMF with a Weighted Variant of Spherical $k$-means}  \label{equiv} 

In this section, we briefly recall how NMF with an additional constraint is equivalent to a fundamental clustering technique (see Equation~\eqref{c1} below): Euclidean $k$-means  \cite{Dhillon07weighted,Ding05onthespectral}. We then observe that relaxing this constraint leads to~\eqref{eq:orth}--\eqref{eq:nonneg}, that is,  ONMF,  which is therefore not exactly equivalent to $k$-means but rather to another problem closely related to spherical $k$-means \cite{Banerjee03}. More precisely, ONMF is equivalent to weighted spherical $k$-means in a particular metric, see Theorem~\ref{Th1}. 
Based on this analysis, we propose a new EM-like algorithm to solve ONMF problems, highlight the differences between $k$-means, spherical $k$-means and ONMF, and illustrate these results on synthetic data sets.

\subsection{Equivalence with Euclidean $k$-means} \label{eqkmeans}

Let $M = (m_{1},\dots,m_{n}) \in \mathbb{R}^{m \times n}_+$ be a nonnegative data matrix whose columns represent a set of $n$ points $\{m_j\}_{j=1}^{n} \in \mathbb{R}^{m}_+$. Solving the clustering problem means finding a set $\{\pi_i\}_{i=1}^{k}$ of $k$ disjoint clusters: 
\[ 
\pi_i \subseteq \{ 1, 2, \dots, n \} \; \forall i, \quad \cup_{1 \leq i \leq k} \pi_i = \{ 1, 2, \dots, n \}, 
\]
\[  
\text{ and }  \quad \pi_i \cap \pi_j = \emptyset \; \forall i \neq j, 
\]
such that each cluster $\pi_i$ contains objects as similar as possible to each other  according to some quantitative criterion. When choosing the Euclidean distance, we obtain the $k$-means problem, which can be formulated as follows \cite{Dhillon07weighted}: 
\begin{equation} \nonumber 
\min_{\{\pi_i\}_{i=1}^{k}} \sum_{i=1}^{k} \sum_{j \in \pi_{i}} \| m_j - c_i \|^{2}, \quad 
\end{equation} 
\[ 
\text{where} \;\; c_i = \frac{\sum_{j \in \pi_i}m_j}{| \pi_i |} \text{ are the cluster centroids}. 
\]
Equivalently, we can define a binary cluster indicator matrix $B \in \{0,1\}^{k \times n}$ as follows:
\[ 
B = \{b_{ij}\}_{k \times n} \quad \text{ where } \quad  b_{ij} = 1 \iff j \in \pi_i. 
\]
Disjointness of clusters $\pi_i$ means that rows of $B$ are orthogonal, i.e., $B B^T$ is diagonal. Therefore we can normalize them to obtain an orthogonal matrix $V = \{v_{ij}\}_{k \times n} = (BB^T)^{-\frac{1}{2}}B$ (a weighted cluster indicator matrix) which satisfies the following condition: 
\begin{align}  
& \text{There exists a set of clusters $\{\pi_i\}_{i=1}^{k}$ such that }   \nonumber \\
& v_{ij} =  
\begin{cases}
    \frac{1}{\sqrt{|\pi_i|}}, & \text{if $j \in \pi_{i}$},  \\
    0, & \text{otherwise}. 
  \end{cases}  \label{c1} \tag{c1} 
\end{align} 
It has been shown in \cite{Ding05onthespectral} that the NMF problem with matrix $V$ satisfying condition  \eqref{c1}: 
\begin{equation} \label{eq:obNMF} 
\min_{U \geq 0, V \geq 0} \| M - UV \|_F^{2} \; \text{ s.t. } \quad
\text{ $V$ satisfies \eqref{c1} }, 
\end{equation}
 is equivalent to $k$-means. In fact, since $V$ in problem \eqref{eq:obNMF} is a normalized indicator matrix which satisfies $v_{ij} = |\pi_i|^{-\frac{1}{2}} \iff j \in \pi_i$, we have 
\begin{eqnarray*} 
\| M - UV \|_F^{2} & = & \sum_{j=1}^n || m_j - \sum_{i=1}^{k} u_i v_{ij} ||^{2}  \\
& = & \sum_{i=1}^{k} \sum_{j \in \pi_i} \| m_j - u_i v_{ij} \|^{2} \\
& =  & \sum_{i=1}^{k} \sum_{j \in \pi_{i}} \| m_j - u_i \frac{1}{\sqrt{|\pi_i|}} \|^{2}, 
\end{eqnarray*} 
which implies that, at optimality, each column $u_i$ of $U$ must correspond (up to a multiplicative factor) to a cluster centroid with $u_i = \sqrt{|\pi_i|} \, c_{i} = \frac{\sum_{j \in \pi_i}m_j}{\sqrt{| \pi_i |}}$ $\forall i = 1,\dots,k$.

\subsection{ONMF and a Weighted Variant of Spherical $k$-means}

Let us now define a condition weaker than \eqref{c1}:
\[ \label{c2} \tag{c2}
VV^T = I_k \quad \text{ and } \quad V \geq 0.
\]
It can be easily checked that $\eqref{c1} \Rightarrow \eqref{c2}$ while $\eqref{c2} \nRightarrow \eqref{c1}$. 
The difference between conditions \eqref{c1} and \eqref{c2} is that condition \eqref{c2} does not require the rows of $V$ to have their nonzero entries equal to each other. Now, if we only impose  the weaker condition \eqref{c2} on NMF, we obtain a relaxed version of \eqref{eq:obNMF} which, by definition, corresponds to orthogonal NMF: 
\begin{equation} \label{eq:JoNMF}
\min_{U \geq 0, V \geq 0} \| M - UV \|_F^{2} \quad \text{ such that } \quad VV^T= I_k. 
\end{equation} 
In the following, we show the equivalence of problem \eqref{eq:JoNMF} with a particular weighted variant of the spherical $k$-means problem: 

\begin{theorem} \label{Th1} 
For a nonnegative data matrix $M \in \mathbb{R}^{m \times n}_+$, the ONMF problem \eqref{eq:JoNMF} is equivalent to the following weighted variant of spherical $k$-means 
\begin{equation}
\max_{\{\pi_i, u_i \in \mathbb{R}^{m}_+, ||u_i||_2 = 1\}_{i=1}^{k}}  \sum_{i=1}^{k}\sum_{j \in \pi_i} \|m_j\|^2 \left( \frac{m_j^T}{\|m_j\|}u_i\right)^2 , 
\label{eq:JoNMFcos}
\end{equation} 
where $\{\pi_i\}_{i=1}^{k}$ is a set of disjoint clusters. 
\end{theorem} 
\begin{proof}
The claim is that~\eqref{eq:JoNMF} and~\eqref{eq:JoNMFcos} are \emph{equivalent}, i.e., a solution of~\eqref{eq:JoNMF} is obtained from a solution of~\eqref{eq:JoNMFcos} by means of elementary arithmetic operations, and vice-versa.

First, without loss of generality, we assume that $k$ is sufficiently small so that the solutions $U$ of~\eqref{eq:JoNMF} do not have vanishing columns. We then redefine ``$U\geq0$'' (resp.\ ``$V\geq0$'') to mean that $U$ (resp.\ $V$) is nonnegative without vanishing columns (resp.\ rows). This redefinition does not alter the the solutions of~\eqref{eq:JoNMF}.

Observe that~\eqref{eq:JoNMF} is equivalent to the following problem obtained by imposing the unit-norm constraint on the columns $\{u_i\}_{i=1}^k$ of $U$ instead of the rows of $V$:
\begin{align} \label{eq:JoNMFd} 
\min_{{U} \geq 0, {V} \geq 0} \| M - {U}{V} \|_F^{2}  \; \text{ s.t. } & \; ({V}{V}^T)_{ij} = 0 \; \forall i \neq j \; \text{ and } \; \|{u}_i\|  = 1 \, \forall i. 
\end{align} 
Indeed, since the function $\psi:(U,V)\mapsto (UD^{-1},DV)$ with $D=\text{diag}(\|u_1\|,\dots,\|u_k\|)$ is a homeomorphism from the feasible set of~\eqref{eq:JoNMF} onto the feasible of~\eqref{eq:JoNMFd} that does not modify the objective value, it is readily seen that if $(U,V)$ is a solution of~\eqref{eq:JoNMF}, then $\psi(U,V)$ is a solution of~\eqref{eq:JoNMFd}; proving the reverse direction is equally straightforward.

It remains to show equivalence between~\eqref{eq:JoNMFd} and~\eqref{eq:JoNMFcos}. 
We say that a partition $\pi = \{\pi_i\}_{i=1}^{k}$ and a matrix $V$ (with $k$ rows) are \emph{compatible}, which we write $V\sim\pi$, if the inclusion $j\in\pi_i$ holds whenever $V_{ij}\neq0$. Using this notion, we first notice the crucial fact that $V\sim\pi$ for some $\pi$ if and only if each column of $V$ has at most one nonzero element (at position $(i,j)$ where $i$ is determined by $\pi_i\ni j$). Defining $\mathcal{V}$ to be the feasible set for $V$ in~\eqref{eq:JoNMFd}, it is now clear that $V\sim\pi$ and $V\geq0$ imply that $V\in\mathcal{V}$ and that, in the reverse direction, one can easily check that $V\in\mathcal{V}$ implies the existence of a partition $\pi$ such that $V\sim\pi$.

We can now show that the following four propositions are equivalent, from which the main claim follows.

\begin{enumerate} 
\item \label{it:1} $U$ and $V$ minimize $\|M-UV\|_F^2$ subject to $U\geq0$, $\|u_i\|=1 \, \forall i$, $V\geq0$, $(VV^T)_{ij} = 0 \, \forall i \neq j$, 
\\ and $\pi$ is obtained by elementary operations to satisfy $V\sim\pi$.
\item \label{it:2} $U$, $V$, and $\pi$ minimize $\|M-UV\|_F^2$ subject to $U\geq0$, $\|u_i\|=1 \, \forall i$, $V\geq0$, $V\sim\pi$.
\item \label{it:3} $U$, $V$, and $\pi$ minimize $\sum_{i=1}^{k} \sum_{j \in \pi_i} \| m_j - u_i v_{ij} \|^{2}$ subject to $U\geq0$, $\|u_i\|=1 \, \forall i$, $V\geq0$, $V\sim\pi$.
\item \label{it:4} $U$ and $\pi$ maximize $\sum_{i=1}^{k} \sum_{j \in \pi_i} \left( m_j^Tu_i\right)^2$ subject to $U\geq0$, $\|u_i\|=1 \, \forall i$,  
\\ and $V$ is obtained by the elementary operations
\[
\begin{cases}
v_{ij}=0 & \text{if $j\notin\pi_i$}
\\ v_{ij}=m_j^Tu_i & \text{if $j\in\pi_i$}.
\end{cases}
\] 
\end{enumerate} 
The equivalence between~\ref{it:1} and~\ref{it:2} follows from the discussion in the previous paragraph. The equivalence between~\ref{it:2} and~\ref{it:3} follows from a rewriting of the objective function made possible by the constraints. 

Finally, the equivalence between~\ref{it:3} and~\ref{it:4} is established as follows. Referring to~\ref{it:3}, given a feasible $U$ and $\pi$, we have for each term $\|m_j - u_i v_{ij} \|^{2}$ that the optimal $v_{ij}^*$ is given by: 
\begin{eqnarray} \label{eq:optCoeffs}
v_{ij}^* 
& = & 
\argmin_{x \geq 0} \|m_j -  u_i x \|^2 \\
 & = &
  \argmin_{x \geq 0} \left( m_j^Tm_j - 2x m_j^Tu_i + x^2\right)  \nonumber \\
  & = & 
  m_j^Tu_i, \; 1 \leq i \leq k, j \in \pi_i. \nonumber  
\end{eqnarray}
(Observe that $m_j^Tu_i\geq0$ in view of the nonnegativity of $M$ and $U$.) 
Backsubstituting the optimal coefficients~\eqref{eq:optCoeffs} in~\ref{it:3}, we have that $U$ and $\pi$ of~\ref{it:3} minimize the function
\begin{align}
& \sum_{i=1}^{k}\sum_{j \in \pi_i}\|m_j - \left( m_j^Tu_i\right)u_i\|^{2}   \notag\\  
 =& 
\sum_{i=1}^{k}\sum_{j \in \pi_i}  \Big( m_j^Tm_j - 2\left( m_j^Tu_i\right)^2 + \left( m_j^Tu_i\right)^2 \Big)  \notag \\
=& \sum_{i=1}^{k}\sum_{j \in \pi_i}  - \left( m_j^Tu_i\right)^2 + \text{cst} .   \notag
\end{align} 
Hence they maximize the function 
\begin{equation}
\sum_{i=1}^{k}\sum_{j \in \pi_i} \left( m_j^Tu_i\right)^2. \label{eq:JoNMFwcos}
\end{equation} 
This shows that~\ref{it:3} implies~\ref{it:4}, and the converse is readily established by contradiction. 
\end{proof}

It is insightful to compare formulation \eqref{eq:JoNMFcos} of ONMF with the spherical $k$-means problem \cite{Banerjee03}, which is a variant of $k$-means where both data points and centroids are constrained to have unit norm: 
\begin{align} 
& \min_{\{\pi_i, u_i\}_{i=1}^{k}} \sum_{i=1}^{k}\sum_{j \in \pi_i} \left\|\frac{m_j}{||m_j||} -  u_i \right\|^2 \; \; \text{ s.t. } \; \; 
||u_i|| = 1,  \notag\\ 
 & \equiv \quad \max_{\{\pi_i, u_i\}_{i=1}^{k}}  \sum_{i=1}^{k}\sum_{j \in \pi_i} \frac{m_j^T}{||m_j||} u_i  \; \; \text{ s.t. } \; \; 
||u_i|| = 1\label{eq:Jspkm}.
\end{align}
Note that both problems \eqref{eq:JoNMFcos} and \eqref{eq:Jspkm} relate to maximizing the cosines of the angles between $u_i$ and the data points from the corresponding cluster. However, we observe that: 

\begin{itemize} 

	\item 
	
Because of coefficients $||m_i||^2$, problem \eqref{eq:JoNMFcos} is sensitive to the norm of the data points,  as opposed to  spherical $k$-means \eqref{eq:Jspkm} which only depends on their direction;  
	
	\item Even for normalized data points (i.e., $||m_i|| = 1$ $\forall i$), problem \eqref{eq:JoNMFcos} is similar but not equivalent to spherical $k$-means \eqref{eq:Jspkm} because it tries to maximize \emph{the sum of squares} of the cosines (instead of their sum). 

\item Contrarily to problem~\eqref{eq:JoNMFcos}, spherical $k$-means \eqref{eq:Jspkm} does not require nonnegativity of $u_i$'s, although it will clearly hold at optimality when data points $m_j$ are nonnegative.
\end{itemize}


\subsection{Which Model Should Be Used: $k$-means, spherical $k$-means or ONMF?}  \label{which}

\begin{figure*}[ht] 
\centering
\subfigure[Well-separated clusters]{\includegraphics[width=6.5cm]{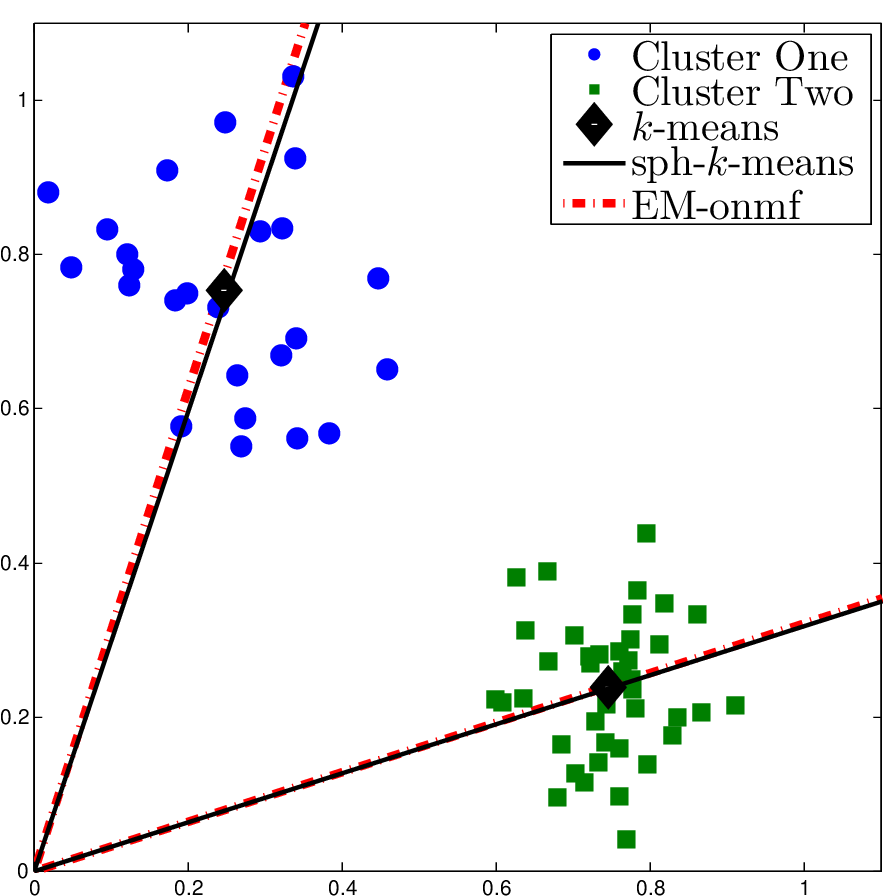}} 
\hspace{2cm}
\subfigure[In-line clusters]{\includegraphics[width=6.5cm]{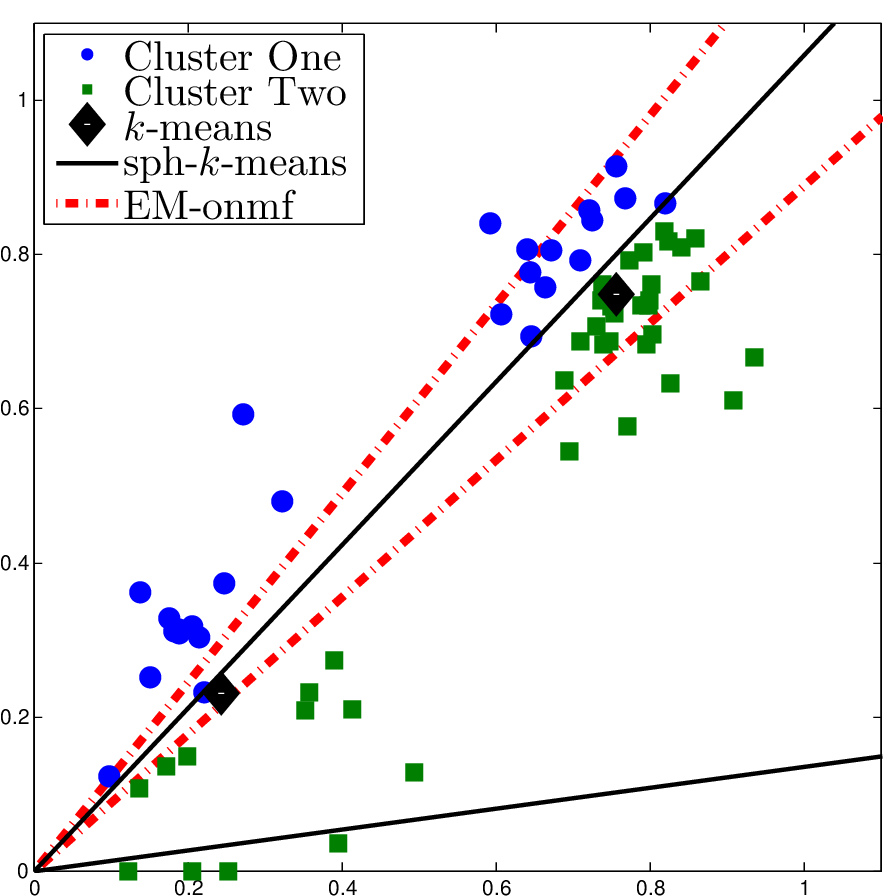}}
\caption{Comparison of $k$-means, standard spherical $k$-means and ONMF. Diamonds are  cluster centroids found by $k$-means, continuous lines are spherical $k$-means centroid directions while dashed lines are ONMF centroid directions. 
Circles and squares are data points as clustered by ONMF. 
As expected, $k$-means is not sensitive to the alignment of the clusters as opposed to spherical $k$-means and ONMF. 
 On the left figure (a), the clusters are well separated and the three techniques perform similarly. On the right figure (b),  the directional effect is clearly visible for both ONMF and spherical $k$-means. However, there is an important difference between the two: ONMF is more sensitive to the data points with larger norm, 
 while 
 spherical $k$-means treats all the points the same way 
 (including the ones from the  lower left cluster with smaller norm but wider angular distribution)  
and its centroids are therefore further apart from each other.} 
\label{compa}
\end{figure*}

Based on the analysis of ONMF from Section~\ref{equiv} (in particular, Theorem~\ref{Th1}), we explain in this section in which situations ONMF should be preferred to $k$-means and spherical $k$-means. 
This issue can be settled by addressing the following two questions: 


\begin{enumerate} 

\item \emph{Should scaling of the data points influence the cluster assignment?} \vspace{0.1cm}

Given the cluster centroids, 
spherical $k$-means and ONMF are invariant to scaling in the sense that, for any $\alpha > 0$, a data point $x$ and its scaling $\alpha x$ will be assigned to the same cluster (the one minimizing the angle between $x$ and the cluster centroid; see Section~\ref{EMalgo}). On the contrary, $k$-means is very sensitive to scaling as it assigns data points to clusters based on  distances; see Figure~\ref{compa} for an illustration. 
In practice, there are many situations where cluster assignment should be independent of scaling so that spherical $k$-means and ONMF should be preferred to $k$-means. For example, in document classification, two documents discussing the same topic will roughly be multiple of one another (scaling depending then on the relative lengths of the documents), and, in hyperspectral imaging, pixels containing the same material will have their spectral signatures multiple of one another (scaling depending on the relative illumination conditions; see Section~\ref{nr} for numerical experiments).  \vspace{0.2cm}

\item \emph{Does the noise added to each data point depend on its norm?} \vspace{0.1cm}

Spherical $k$-means is invariant under normalization of the data points (see Equation~\eqref{eq:Jspkm}) while ONMF gives more importance to data points with larger norm. For example, if the noise added to each data point is independent of its norm, ONMF should be preferred. 
In fact, in that situation, data points with larger norm are relatively contaminated with less noise hence should be given more importance. 
Another example is when data points with larger norm are statistically more significant. This is usually the case for example in document classification: assuming that each document discusses only one topic and that each topic is a distribution over the words, longer documents represent a larger sample of the corresponding topic distribution and should be given more importance. (In fact, most document classification software typically discards very short documents: ONMF implicitly takes care of this issue by giving less importance to shorter documents.)  In hyperspectral imaging, background pixels contain mostly noise and should not be given too much importance; hence ONMF should also be preferred in this situation.  


\end{enumerate} 
Figure 1 displays a comparison between $k$-means, spherical $k$-means and ONMF on two simple examples.  
 \\

\begin{figure*}[ht] 
\centering
\includegraphics[width=8cm]{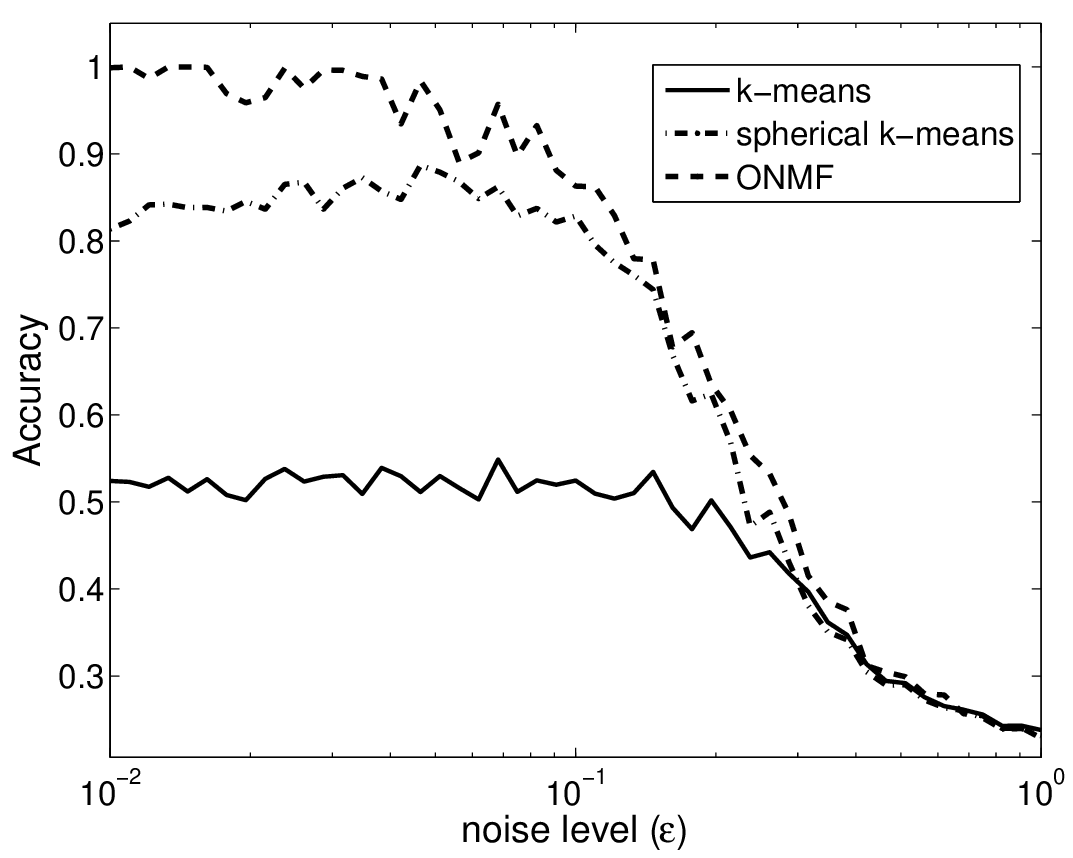}
\caption{Comparison of $k$-means, standard spherical $k$-means and ONMF on synthetic data sets.} 
\label{syntdat} 
\end{figure*}

To conclude, ONMF should be preferred to both $k$-means and spherical $k$-means when 
\begin{itemize}
\item[(1)] Scaling should not affect cluster assignment, and 
\item[(2)] Data points with larger norm are more reliable and should be given more importance. 
\end{itemize}
To illustrate this, we generate several synthetic data sets as follows. Each data sets has six clusters $\{\pi_i\}_{i=1}^{6}$, each containing $100-(i-1)10$ data points for a total of $450$ data points. Each cluster centroid $u_i \in \mathbb{R}^{10}$ $1 \leq i \leq 6$ is generated uniformly at random in the unit cube $[0,1]^6$. 
Each data point $m_j$ $1 \leq j \leq 450$ is a multiple of its corresponding cluster centroid: $m_j = \alpha u_k$ where $\alpha > 0$ is picked uniformly at random in the interval $[0.1,1]$. 
Hence, because of this scaling, $k$-means is not appropriate and will perform poorly. 
Each data point is then perturbed by some additive noise with fixed magnitude (i.e. independent of the norm of the data point, which should lead to ONMF performing better than spherical $k$-means). Concretely, each noise entry is drawn from a normal distribution with zero mean and fixed standard deviation $\epsilon$. Finally, the negative entries of each data point are set to zero to obtain a nonnegative input matrix (note that this can only reduce the noise). 
Letting $\{\pi'_i\}_{i=1}^{6}$ be the clusters extracted by an algorithm and $\{\pi_i\}_{i=1}^{6}$  be the true clusters, the accuracy is defined as 
\begin{equation} \label{accu}
\text{Accuracy} 
= 
\max_{P \in [1,2,\dots,k]} \frac{1}{450} \left( \sum_{i=1}^6  | \pi_i \cap \pi'_{P(i)} | \right)  \; \in \;  [0,1], 
\end{equation}
where $[1,2,\dots,k]$ is the set of permutations of $\{1,2,\dots,k\}$. 
For each noise level $\epsilon \in [0,1]$, we generate ten synthetic data sets as described above, and Figure~\ref{syntdat} reports the average accuracy of each algorithm using ten random initializations (except for ONMF which was solved using ONP-MF, which does not use random initialization, see Section~\ref{npomf}). 
As expected, we observe that ONMF outperforms $k$-means and spherical $k$-means. We will use the same synthetic data sets to compare the different ONMF algorithms in Section~\ref{nr}.

\subsection{EM-like Algorithm for ONMF} \label{EMalgo}

We present here a simple EM-like alternating algorithm designed to tackle the ONMF problem (\ref{eq:JoNMF}) based on its equivalence with the weighted variant of spherical $k$-means \eqref{eq:JoNMFcos}. It is very similar to the standard spherical $k$-means algorithm \cite{Banerjee03}, except for the computation of cluster centroids.  Specifically, it starts with an initial set of centroids, either randomly chosen or supplied as initial values. It then alternates between two steps: 
\begin{enumerate}
\item Given cluster centroids $\{u_i\}_{i=1}^k$, choose $\{\pi_i\}_{i=1}^k$ assigning each point to its closest cluster: 
\begin{align*}
j \in \pi_i \Rightarrow i & \in \argmax_{1 \leq \ell \leq k} {\left(m_{j}^Tu_{\ell}\right)^2} \\
& = \argmax_{1 \leq \ell \leq k} {\left(m_{j}^Tu_{\ell}\right)}.
\end{align*}
Notice that this step is exactly equivalent to the one of standard spherical $k$-means \cite{Banerjee03}. 
\item Given the clustering $\{\pi_i\}_{i=1}^k$, compute the new optimal cluster centroids $\{u_i\}_{i=1}^k$ as follows. Define matrix $M_i \in \mathbb{R}^{m \times |\pi_i|}$ as the submatrix of $M$ containing the columns belonging to cluster $\pi_i$. 
We have to solve problem \eqref{eq:JoNMFwcos} with respect to the $u_i$'s: 
\[
\max_{\{u_i \geq 0, ||u_i|| = 1\}_{i=1}^{k}}  \sum_{i=1}^{k}\sum_{j \in \pi_i} \left( m_j^Tu_i\right)^2 = \sum_{i=1}^{k} ||M_i^Tu_i||_2^2 . 
\]
There are $k$ independent problems: each $u_i$ must maximize the term $||M_i^Tu_i||_2^2$. The optimal solution $u_i^*$ is given by the dominant left singular vector of $M_i$ associated with $\sigma_1(M_i)$, the largest singular value of $M_i$:
\[
u_i^* = \argmax_{||u||_2 = 1} ||M_i^Tu||_2^2 = \argmax_{||u||_2 = 1}  u^T M_i M_i^T u,  
\]
for which we have $||M_i^Tu_i^*||_2 = \sigma_1(M_i) = ||M_i||_2$. Moreover, since $M_i \ge 0$, the Perron-Frobenius theorem guarantees that $u_i^*$ can always be chosen to be nonnegative.
\end{enumerate}
Algorithm~\ref{algomod}, referred to as EM-ONMF, implements this procedure.  We will see in the last section that, despite its simplicity, it works well for text clustering tasks.   
\begin{algorithm} \label{algomod}
\SetKwInOut{Input}{input}\SetKwInOut{Output}{output}
\Input{ Nonnegative data matrix $M$, and initial centroids $\{u_i\}_{i=1}^{k}$.}
\Output{  Clustering of the points $\{\pi_i\}_{i=1}^{k}$, with the corresponding centroid directions $\{u_i\}_{i=1}^{k}$.}
\BlankLine
\While{not converged}{
$\{\pi_i\}_{i=1}^{k} \leftarrow \emptyset$\;
\For{$j\leftarrow 1$ \KwTo $n$}{
find  $i \in \argmax_{1 \leq \ell \leq k} {\left(m_j^T u_{\ell}\right)}$ and update cluster $\pi_{i} = \pi_{i} \cup \{ j \}$. \;
}
\textbf{if} $\pi_i = \emptyset$ for some $i$ \textbf{then} randomly transfer a point to cluster $\pi_i$ 
\textbf{end} \;
\For{$i\leftarrow 1$ \KwTo $k$}{
$u_i$ $\leftarrow$ (any) nonnegative dominant singular vector of the data submatrix $M_i = M(:,\pi_i)$. \;
}
}
\caption{EM-like Algorithm for ONMF (EM-ONMF)} 
\end{algorithm} 
Note that Algorithm~\ref{algomod} does not explicitly provide a solution to ONMF. However, a candidate solution of~\eqref{eq:JoNMFd} can be obtained by taking $U=\begin{bmatrix}u_1 & \dots & u_k\end{bmatrix}$ and selecting $V$ according to~\eqref{eq:optCoeffs}, from which a candidate solution of ONMF~\eqref{eq:JoNMF} is readily obtained as $(UD,D^{-1}V)$ where $D=\mathrm{diag}(\|V(1,:)\|_2,\dots,\|V(k,:)\|_2)$. \\

It is interesting to relate this with the original ONMF problem \eqref{eq:JoNMFd}: given a partitioning $\{\pi_i\}_{i=1}^{k}$, let us denote $w_i = (v_{ij})_{j \in \pi_i}$ the subvector containing only the positive entries of the $i^{\text{th}}$ row of $V$. Then, 
\begin{align*}
\| M - {U}{V} \|_F^{2}  &=   \sum_{i=1}^{k} \sum_{j \in \pi_{i}} \| m_j - u_i v_{ij} \|^{2} \\
											  &=   \sum_{i=1}^{k} ||M_i - u_i w_i^T||_F^2, 			
\end{align*}
so that the optimal ($u_i$,$w_i$) must be an optimal solution of 
\begin{equation} \label{probMi}
\min_{||u_i|| = 1, u_i  \geq 0, w_i \geq 0}  ||M_i - u_i w_i^T||_F^2. 
\end{equation}
Each of these problems looks for the best nonnegative rank-one approximation of a nonnegative matrix (i.e., a rank-one NMF problem). This in turn can be solved 
by combining the Eckart-Young and Perron-Frobenius theorems: taking 
 the first rank-one factor generated by the singular value decomposition (SVD) (making sure it is nonnegative in case of non-uniqueness) leads to a minimum value for \eqref{probMi} equal to $||M_i||_F^2 - \sigma^2_1(M_i)$. Therefore, solving ONMF amounts to finding a partitioning $\{\pi_i\}_{i=1}^{k}$ such that the sum of squares of the first singular values of submatrices $M_i$'s is maximized, that is,  
 \emph{ONMF problem \eqref{eq:JoNMF} is equivalent to 
$\max_{\{\pi_i\}_{i=1}^{k}} \sum_{i=1}^k \sigma_1^2(M_i)$. } 

\section{Augmented Lagrangian Method for ONMF}  \label{npomf} 

In this section, we present an alternative approach to solve ONMF problems\footnote{Our code for the proposed algorithm is available at \url{https://bitbucket.org/filp/onmf/src}.}. Typically, ONMF algorithms strictly enforce nonnegativity for each iterate while trying to achieve orthogonality at the limit. This can be done using a proper penalization term \cite{Ding06orthogonalnonnegative}, a projection matrix formulation \cite{Yang10linearand} or by choosing a suitable search direction \cite{Choi08algorithmsfor}. 
We propose here a method working the opposite way: at each iteration, a (continuous) projected gradient scheme is used to ensure that the $V$ iterates are 
 orthogonal (but not necessarily nonnegative). 

\begin{figure*}[ht] 
\centering
\begin{tabular}{ccc}
\includegraphics[width=7cm]{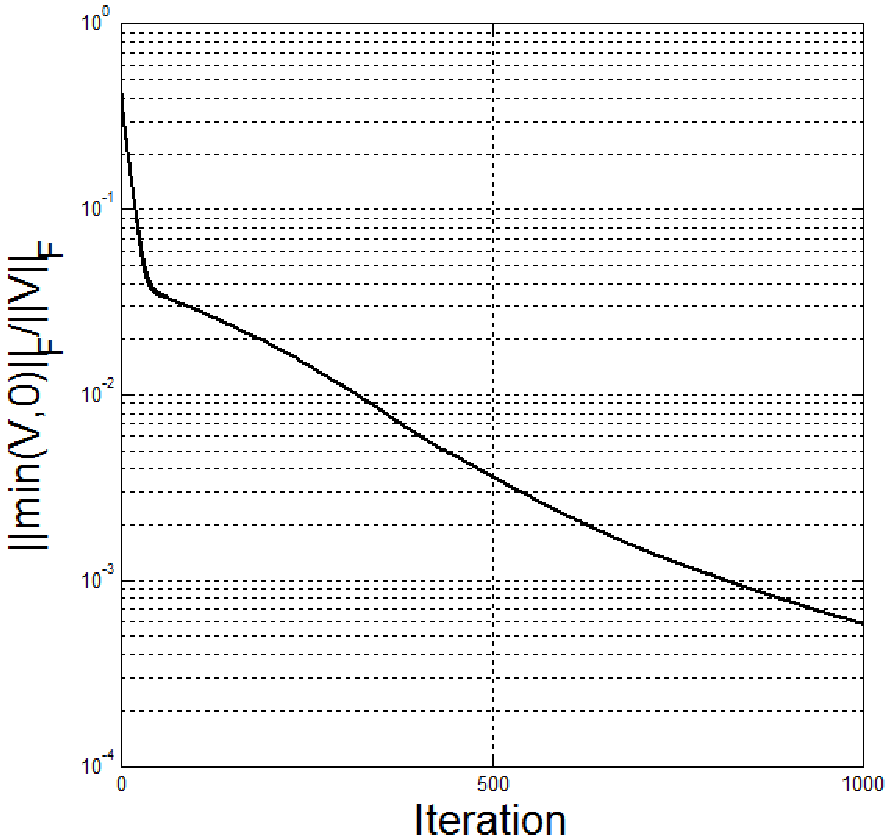}& \qquad \qquad \qquad & \includegraphics[width=7cm]{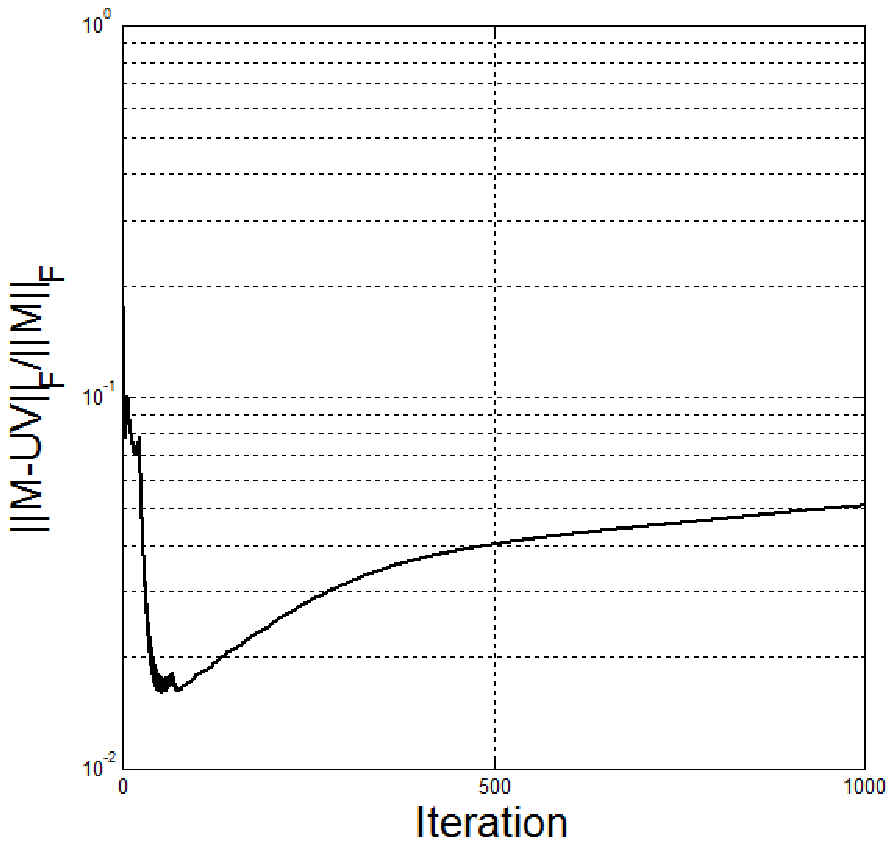} \\
\end{tabular}
\caption{Convergence of Alg.~\ref{ONPMF} for the Hubble data set (left: constraint residual, right: approximation error).}
\label{conv}
\end{figure*}

Nonnegativity constraints in the ONMF formulation \eqref{eq:JoNMF} will be handled using the following augmented Lagrangian, defined for a matrix of Lagrange multipliers $\Lambda \in \mathbb{R}^{k \times n}_+$  associated to the nonnegativity constraints:
\begin{equation} \label{lag}
L_{\rho}(U,V,\Lambda)  = \frac{1}{2}||M-UV||_F^2 + \langle \Lambda,-V \rangle
 + \frac{\rho}{2}|| \min(V,0) ||_F^2, 
\end{equation}
where $\rho$ is the quadratic penalty parameter. 
Ideally, we would like to solve the Lagrangian dual 
\[
\max_{\Lambda \geq 0} f(\Lambda) \quad \text{ where } \quad f(\Lambda) = \min_{U \geq 0, VV^T = I_{k}} L_{\rho}(U,V,\Lambda). 
\]
Observe that,
regardless of the value of $\rho$, the
solutions $(U,V)$ of the ONMF problem~\eqref{eq:ONMF} are the solutions of
\[
\min_{U\geq0, VV^T=I_k} \max_{\Lambda\geq0} L_{\rho}(U,V,\Lambda).
\] 
We propose here a simple alternating scheme to update variables $U$, $V$, $\Lambda$ while, as announced, explicitly enforcing $U\geq0$ and $VV^T=I_k$:
\begin{enumerate}

\item For $V$ and $\Lambda$ fixed, the optimal $U$ can be computed by solving a nonnegative least squares problem 
$
U \leftarrow \argmin_{X \in \mathbb{R}^{m \times k}_+} \|M - XV\|^{2}_F. 
$
We use the efficient active-set method  
proposed in\footnote{Available at \url{http://www.cc.gatech.edu/~hpark/}.}  \cite{KP11}. 

\item For $U$ and $\Lambda$ fixed, we update  matrix $V$ by means of a projected gradient scheme. 
Computing the projection of a matrix $\hat{V}$ onto the feasible set of orthogonal matrices, known as the Stiefel manifold\footnote{The Stiefel manifold is the set of all $n \times k$ orthogonal matrices, i.e., St($k$,$n$) $= \{ X \in \mathbb{R}^{n \times k} : X^TX = I_k$\}.}, amounts to solving the following problem: 
\[
\text{Proj}_{St}(\hat{V}) = \argmin_{X} || \hat{V} - X ||_F^2 \quad \text{ such that } \; XX^T = I_{k}, 
\]
whose optimal solution $X^*$ can be computed in closed form from the unitary factor of a polar decomposition of $\hat{V}$, see, e.g., \cite{HJ90, AbsMal2012}. Our projected gradient scheme then reads: 
\[
V \leftarrow \text{Proj}_{St} \Big(  V  - \beta \nabla_{V} L_{\rho}(U,V,\Lambda) \Big),  
\]
where the step length $\beta$ is chosen with a backtracking line search similar to that in \cite{L07} (step length is increased as long as there is a decrease in the objective function, and decreased otherwise).

\item Finally, Lagrange multipliers are updated in order to penalize the negative values of $V$:
\[
\Lambda \leftarrow  \max\left({0}, \Lambda - \alpha V \right). 
\]
where $\alpha$ is the step length. 
As $-V$ is the gradient of function $\Lambda\mapsto L_{\rho}( {U},{V},\Lambda)$, this update is a (projected) gradient step with step length $\alpha$. We choose a predefined step length sequence $\alpha = \alpha_{0}/t$, where $t$ is the iteration counter and $\alpha_{0} > 0$ is a constant parameter, that satisfies the usual ``square summable but not summable'' condition of online gradient methods~\cite[(5.1)]{Bot1998}.  
\end{enumerate}
To initialize the algorithm, we set $\Lambda$ to zero and choose for the columns of $V$ the first $k$ right singular vectors of the data matrix $M$ (which can be obtained with SVD)\footnote{To overcome the sign ambiguity of each row of $V^{(0)}$, we flip its sign if the $\ell_2$-norm of its negative entries is larger than the $\ell_2$-norm of its positive entries.}. 
Quadratic penalty parameter $\rho$ is initially fixed to a small value $\rho_0$ and then increased geometrically after each iteration. Alg.~\ref{ONPMF} implements this procedure, which we refer to as Orthogonal Nonnegatively Penalized Matrix Factorization (ONP-MF). 
\begin{algorithm} 
\caption{Orthogonal nonnegatively penalized matrix factorization (ONP-MF)} 
\SetKwInOut{Input}{input}\SetKwInOut{Output}{output}
\Input{ A nonnegative data matrix $M$, the number of clusters $k$, $\alpha_{0} > 0$,  $\rho_0 > 0$ and $C > 1$.} 
\Output{ The centroid matrix $U$, and the cluster assignment matrix $V$.} 
\BlankLine
Initialize $\Lambda^{(0)} = {0}$, the rows of $V^{(0)}$ with the first $k$ right singular vectors of $M$, and $\rho = \rho_0$. 
\;
\For{$t = 1, 2, \dots$}{
Update $U^{(t)}$ with the optimal solution $U^* = \argmin_{U \geq 0} ||M-UV^{(t-1)}||_F^2$.\;
Update $V^{(t)}$ with projected gradient and a line search for step $\beta^{(t)}$:
\[
V^{(t)} \leftarrow \text{Proj}_{St} \left[  V^{(t-1)}  - \beta^{(t)} \nabla_{V}L_\rho\bigl( U^{(t)},V^{(t-1)},\Lambda^{(t-1)} \bigr)  \right]. \;
\]
Update Lagrange multipliers (using an approximate subgradient): 
\[
\Lambda^{(t)} \leftarrow  \max\left({0}, \Lambda^{(t-1)} - \frac{\alpha_{0}}{t} V^{(t)} \right). \;
\]
Update $\rho \leftarrow C \rho$. 
}
\label{ONPMF}
\end{algorithm}
We observed that the term $||\min(V,0)||_F$ decreases linearly to zero (as augmented  Lagrangian methods are expected to, see \cite[Th.~17.2]{NR06}) while $||M-UV||_F$ converges to a fixed value, see Figure~\ref{conv} for an example on the Hubble data set (cf.\@ Section~\ref{hu}). A rigorous convergence proof is a topic for further research. 
In fact, it is difficult to analyze an augmented Lagrangian approach when the subproblems are not solved exactly (in our case, using a single loop of a block coordinate descent method) and, as far as we know, such a proof does not exist in the literature (see, e.g., \cite{NR06}) although it is a very popular method. 

%

\section{Numerical Experiments} \label{nr}

In this section, we report some preliminary numerical experiments showing that ONP-MF (Alg.~\ref{ONPMF}) and EM-ONMF (Alg.~\ref{algomod}) perform competitively with two recently proposed methods for ONMF: CHNMF from Choi \cite{Choi08algorithmsfor} and \ProjNMF \, from Yang and Oja \cite{Yang10linearand} (Euclidean variant). 
It should be noted that because ONP-MF is initialized with SVD, its results are deterministic and obtained with just one execution of the algorithm.  
However, it could be argued that the comparison with the other algorithms is not completely fair as CHNMF and \ProjNMF \, are initialized with randomly generated factors. In order to perform a fairer comparison, we also initialize CHNMF and \ProjNMF \, with an SVD-based initialization \cite{Bou08SVD} (SVD cannot be used directly because its factors are not necessarily nonnegative), which will be denoted CH(SVD) and \ProjNMFSVD \, respectively. 
Finally, we also report results from two standard EM clustering algorithms, namely $k$-means and spherical $k$-means (SKM) (see, e.g., \cite{Banerjee03}). 
We will see that EM-ONMF is quite efficient for text clustering tasks (see Section \ref{tc}) while ONP-MF gives very good results for unsupervised image classification tasks (see Sections \ref{hu} and \ref{ima}). 

Parameters for ONP-MF are chosen as follows: 
 $\alpha_{0} = 100$, $\rho_0 = 0.01$ and $C=1.01$ for all data sets. 
 ONMF algorithms are run until a stopping condition is met (see below), or a maximum of  5000 iterations in case of random initializations (for CHNMF and \ProjNMF) and 20000 iterations for the SVD-based initialization (as done in \cite{Bou08SVD}) was reached. The following stopping condition for CHNMF seems to work well in practice\footnote{It seems that $10^{-7}$ is a good trade-off: for example, using $10^{-8}$ instead leads to much larger computational times without significant improvements in clustering accuracy.}:   
\[
\frac{|\,||M-U^{(t+1)}V^{(t+1)}||_F - ||M-U^{(t)}V^{(t)}||_F|}{||M||_F} < 10^{-7}, 
\]
where $t$ is the iteration counter. For \ProjNMF, we use the stopping criterion suggested by its authors\footnote{\iftoggle{use_ORTH_PNMF_instead_of_PNMF}{Code available at \url{http://users.ics.tkk.fi/rozyang/pnmf/index.html}.}{We used our own implementation of the orthogonal PNMF based on a simple modification of the non-orthogonal variant available at \url{http://users.ics.tkk.fi/rozyang/pnmf/index.html}.}}: 
\[
\frac{||V^{\left(t-1\right)}-V^{\left(t\right)}||_F}{||V^{\left(t-1\right)}||_F} < 10^{-5}. 
\]
For ONP-MF, we check whether the current iterate is `sufficiently' nonnegative, using
\[
\frac{||\min\left(V,0\right)||_F}{||V||_F} < 10^{-3}.
\]
All EM-like algorithms, EM-ONMF included, were run until cluster assignment did not change for two consecutive iterations. The initial centroids  were randomly selected among the data points. 
For each experiment, a number of 30 repetitions was executed in random conditions both for ONMF and EM-like algorithms (except for the synthetic data sets where we only performed 10 as in Section~\ref{which}).  In the image experiments, we will display the best solution obtained, i.e, with the lowest error.
All experiments were run on an Intel$^{\textrm{\textregistered}}$ Core{\texttrademark} i7-2630QM quad core CPU @2.00GHz with 8GB of RAM.

\subsection{Synthetic Data Sets} \label{sds} 

\begin{figure*}[ht] 
\centering
\includegraphics[width=17cm]{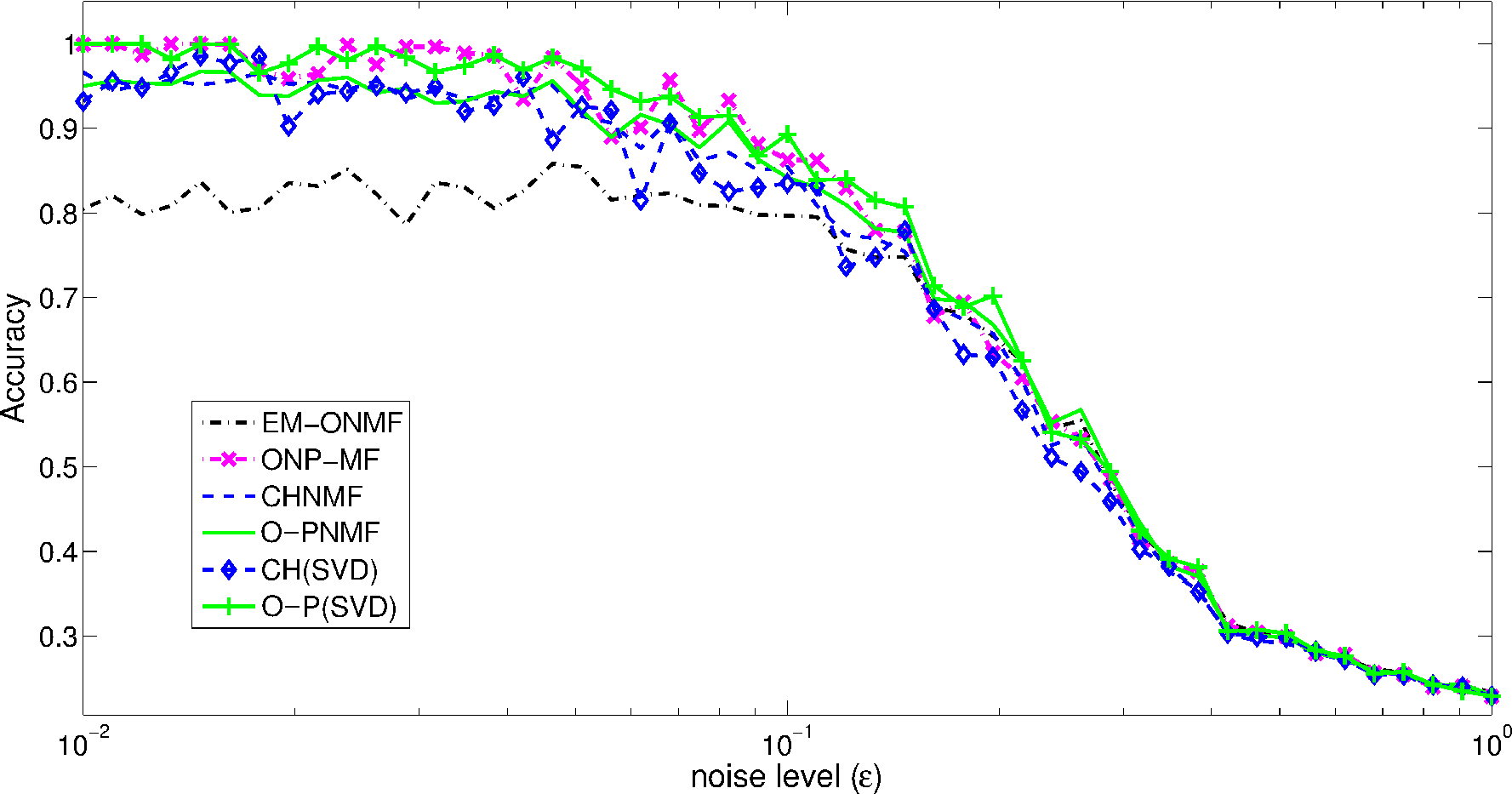}
\caption{Comparison of the different ONMF algorithms on synthetic data sets.} 
\label{syntdat2} 
\end{figure*}

In this section, we perform experiments on synthetic data sets as in Section~\ref{which} in order to compare the different ONMF algorithms. 
 Figure~\ref{syntdat2} reports the average accuracy of each algorithm. 
We observe that 
\begin{itemize}

\item ONP-MF or \ProjNMFSVD \, perform the best among all ONMF algorithms. In particular, they are the only algorithms able to perfectly identify all clusters for small noise levels ($\epsilon = 0.01$). 

\item CHNMF, \ProjNMF \, and CH(SVD) perform similarly, their accuracies being in most cases lower than those of ONP-MF and \ProjNMFSVD. 

\item EM-ONMF performs rather poorly and is in general not able to perform a good clustering (although it is much faster than the other algorithms, see Table~\ref{sdt}). 

\end{itemize}

Table~\ref{sdt} gives for each algorithm the average computational time for one execution on a synthetic data set.
\begin{table}[h!] 
\begin{center}
\caption{Average computational time in seconds for the synthetic data sets.} 
\begin{tabular}{|c|c|c|c|c|c|} 
\hline
  CHNMF & CH(SVD)  & \ProjNMF & \ProjNMFSVD & EM-ONMF &  ONP-MF   \\ 
  2.19  & 3.25     & 2.38 & 2.73   & 0.41    & 3.80 \\ 
\hline
\end{tabular}
\label{sdt} 
\end{center}
\end{table}
ONP-MF is slightly slower but typically obtains one of the best factorizations using only a single deterministic initialization.

\subsection{Text Clustering} \label{tc}

We selected twelve well-known preprocessed document databases described in \cite{ZG05}. Each data set is represented by a term-by-document matrix of varying characteristics, see Table~\ref{dtm}. 
\begin{table}[h!] 
\begin{center}
\caption{Text mining data sets \cite{ZG05}.} 
\begin{tabular}{|c|c|c|c|c|} 
\hline
Data &   $m$  &  $n$ &  r &  $\#\text{nonzero}$   \\ \hline \hline
classic &  7094   & 41681 & 4 & 223839   \\ 
ohscal &   11162  & 11465 & 10  & 674365  \\
hitech &    2301 & 10080 & 6 & 331373 \\  
reviews &   4069  & 18483 & 5 & 758635   \\ 
sports &    8580 & 14870 &  7 & 1091723  \\ 
la1 &    3204 & 31472 &  6 & 484024    \\ 
la2 & 3075 & 31472 & 6 & 455383 \\
k1b & 2340 & 21839 & 6 & 302992 \\
tr11 & 414 & 6429 & 9 & 116613\\
tr23 & 204 & 5832 & 6 &  78609\\
tr41 & 878 & 7454 & 10 & 171509\\
tr45 & 690 & 8261 & 10 & 193605\\
\hline
\end{tabular}
\label{dtm} 
\end{center}
\end{table} 
As a performance indicator, we use the \emph{accuracy}; see Equation~\eqref{accu}. 
We report the average value of the obtained accuracy along with the standard deviations in Table~\ref{tabTxtAcc}. For more than half of the data sets, the average best result was achieved by our algorithms, either EM-ONMF or ONP-MF. Moreover, our algorithms obtain the best performance among ONMF algorithms in ten out of the twelve data sets (being close for the two remaining data sets). 
While EM-ONMF is very fast with a low number of iterations (as the other EM-like algorithms), ONP-MF is in general slower than the other ONMF algorithms (especially CHNMF), and typically requires a larger number of iterations to converge. 

Note that, in this paper, our focus is on the comparison of ONMF algorithms based on the Frobenius norm. 
Comparison with ONMF algorithms using other measures, such as the Kullback-Leibler divergence, and comparison with other topic models (such as LDA \cite{BNJ03}) is a topic for further research. 

\begin{table*}
\begin{center}
\caption{Average accuracy and standard deviation (if applicable) in percent obtained by the different algorithms (in bold, best average performance; underlined, second best).} 
\begin{tabular}{| c || c c | c c c c c c| }
\hline
\hspace{-0.2cm} data set \hspace{-0.2cm} & $k$-means \hspace{-0.2cm} & \hspace{-0.3cm} SKM \hspace{-0.3cm} & \hspace{-0.2cm} CHNMF \hspace{-0.2cm} & \hspace{-0.3cm} CH(SVD) \hspace{-0.3cm} & \hspace{-0.2cm} \ProjNMF \hspace{-0.2cm} 
& \hspace{-0.2cm} \ProjNMFSVD \hspace{-0.2cm} &  EM-ONMF \hspace{-0.2cm} & \hspace{-0.3cm} ONP-MF \hspace{-0.2cm} \\
\hline
\iftoggle{use_ORTH_PNMF_instead_of_PNMF} 
{
classic & \textbf{58.9}   $\pm$6.8 & 56.8   $\pm$4.7  &  55.0  $\pm$1.8  &  55.9 &  50.8 $\pm$1.9   &	  53.9 & \underline{58.8} $\pm$6.9&  53.8\\
ohscal &  28.8   $\pm$3.2 & \textbf{42.8}   $\pm$2.9   &  33.7   $\pm$2.6  &  34.0 &  35.0 $\pm$1.1  &  33.8 & \underline{39.2} $\pm$3.3	&  34.0\\
hitech &  32.2   $\pm$1.6 & \textbf{48.7}   $\pm$3.6&  42.0   $\pm$4.2 &  41.5  &  47.0 $\pm$1.0  & \underline{47.7} & \textbf{48.7} $\pm$5.6 &  47.0\\
reviews & 43.6    $\pm$5.5& \textbf{68.7}   $\pm$6.5 &  52.6   $\pm$9.1 &  49.3 &  55.6 $\pm$6.4  &  52.8 & \underline{63.7} $\pm$10.3 	&  51.0\\
sports &  38.8    $\pm$2.4 &  45.1   $\pm$4.1   & 42.3      $\pm$2.9 & \underline{49.5}  &  44.3 $\pm$3.8  &  49.0 & 	\textbf{50.0} $\pm$6.6	& \textbf{50.0}\\
la1 &     35.0    $\pm$1.9&  48.0   $\pm$4.7   & 53.0      $\pm$5.4 &  44.3 &  55.4 $\pm$5.8  & \underline{60.9} &  50.2 $\pm$7.3 	& \textbf{65.8}\\
la2 &     33.8   $\pm$1.7  &  46.6   $\pm$4.5    & 41.0      $\pm$2.7 &  42.1 &  48.0 $\pm$4.6  & \underline{52.7} &  47.1 $\pm$6.1 	& \textbf{52.8}\\  
k1b &     66.8   $\pm$10.1 &  64.9   $\pm$8.2   & 74.9      $\pm$2.9 & \underline{76.7} &  72.7 $\pm$5.7  &  76.4 &  75.3 $\pm$6.9 	& \textbf{79.0} \\
tr11 & 31.6  $\pm$1.7 &  \textbf{53.0}  $\pm$5.2  &  47.1  $\pm$3.3  &  \underline{50.2}  &  31.5 $\pm$8.5   &  33.8   &  42.4 $\pm$6.3   &  46.1 \\
tr23 &40.8  $\pm$2.5  &  \textbf{42.5}  $\pm$5.2 &  37.0 $\pm$3.6  &  32.8   &  39.2 $\pm$2.0   &  \underline{41.2}   &  40.7 $\pm$4.4   &  40.7 \\
tr41 &41.8 $\pm$6.6  &  \textbf{53.2}   $\pm$5.8  &  46.5 $\pm$5.5 &  42.6   &  35.5  $\pm$7.8  &  43.1   &  \textbf{53.2} $\pm$7.4   &  43.1 \\
tr45 &27.9   $\pm$4.2 &  \textbf{54.2} $\pm$6.2  &  39.2  $\pm$1.7 &  39.3   &  35.7  $\pm$4.7  &  35.1   &  \underline{41.4} $\pm$6.6   &  35.9 \\ 
}
{
classic & \textbf{0.589} & 0.568 & 0.550 & 0.559 & 0.538 &	 0.548 & \underline{0.588}	& 0.538\\
ohscal &  0.288 & \textbf{0.428} & 0.337 & 0.340 & 0.344 & 0.353 & \underline{0.392}	& 0.340\\
hitech &  0.322 & \textbf{0.487} & 0.420 & 0.415  & 0.449 & \underline{0.471} & \textbf{0.487} & 0.470\\
reviews & 0.436 & \textbf{0.687} & 0.526 & 0.493 & 0.519 & 0.503 & \underline{0.637} 	& 0.510\\
sports &  0.388 & 0.451 & 0.423 & \underline{0.495}  & 0.430 & 0.410 & 	\textbf{0.500}	& \textbf{0.500}\\
la1 &     0.350 & 0.480 & 0.530 & 0.443 & 0.549 & \textbf{0.660} & 0.502	& \underline{0.658}\\
la2 &     0.338 & 0.466 & 0.410 & 0.421 & 0.470 & \underline{0.510} & 0.471	& \textbf{0.528}\\  
k1b &     0.668 & 0.649 & 0.749 & \underline{0.767} & 0.746 & 0.757 & 0.753	& \textbf{0.790} \\
}
\hline
\end{tabular} 
\label{tabTxtAcc}
\end{center}
\end{table*}


\subsection{Hyperspectral Unmixing} \label{hu}

A hyperspectral image is a set of images of the same object or scene taken at different wavelengths. Each image is acquired by measuring the reflectance (i.e., the fraction of incident electromagnetic power reflected) of each individual pixel at a given wavelength. The aim is to classify the pixels in different clusters, each representing a different material. 
We want to cluster the columns of a wavelength-by-pixel reflectance matrix so that each cluster (a set of pixels) corresponds to a particular type of material.

\subsubsection{Hubble Telescope} 

\begin{figure*}[ht] 
\centering
\begin{tabular}{c}
\begin{tabular}{ccc}
\includegraphics[width=2.80cm]{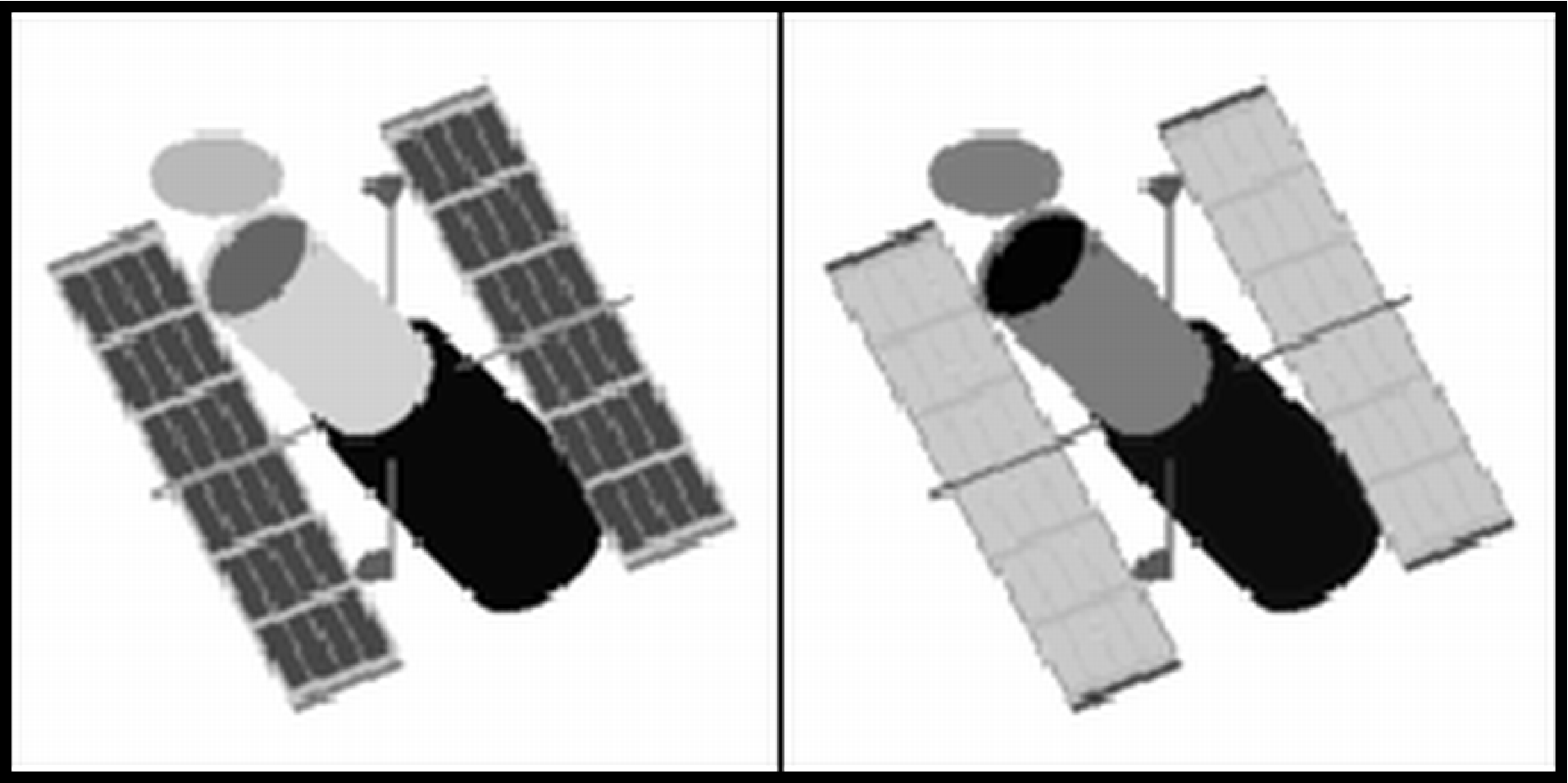} & \qquad & \includegraphics[width=11cm]{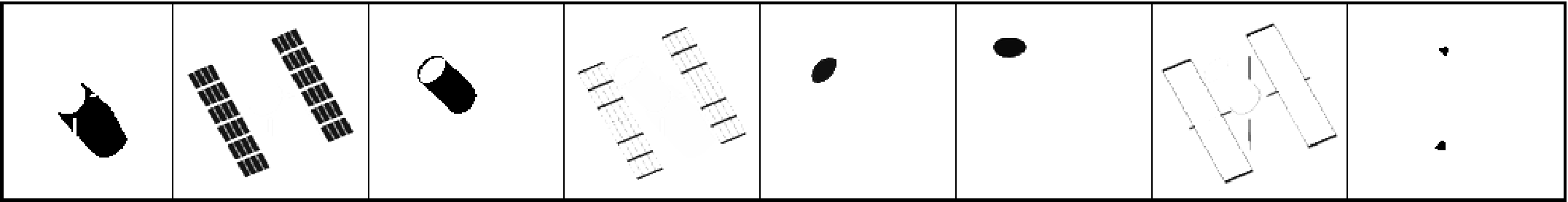}\\
\end{tabular}\\
\includegraphics[width=11cm]{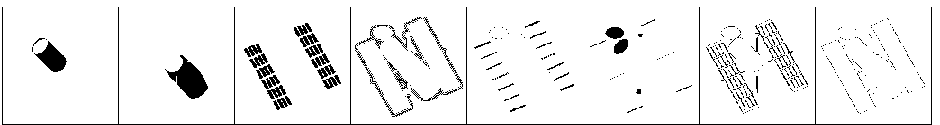} \\
\includegraphics[width=11cm]{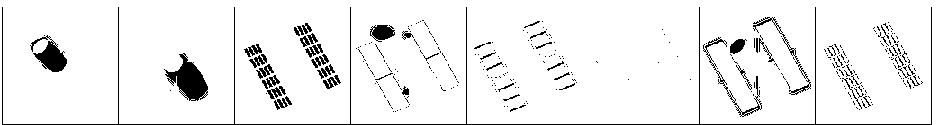} \\
\includegraphics[width=11cm]{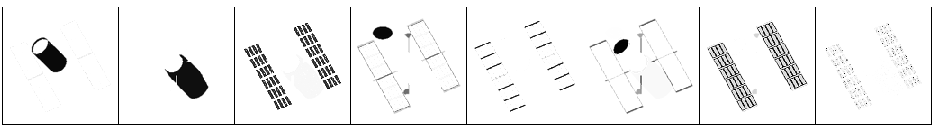} \\
\includegraphics[width=11cm]{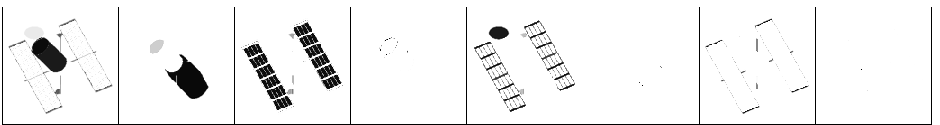} \\
\iftoggle{use_ORTH_PNMF_instead_of_PNMF}
{
\includegraphics[width=11cm]{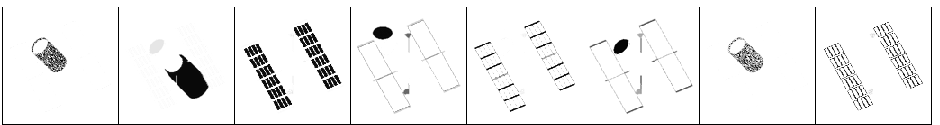} \\
\includegraphics[width=11cm]{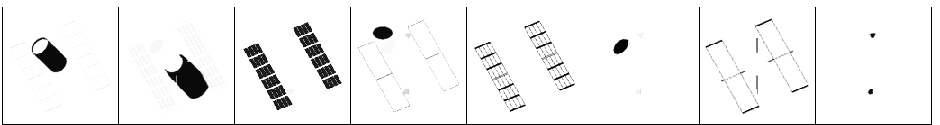} \\
}
{
\includegraphics[width=11cm]{HUBBLE_PNMFRAND.eps} \\
\includegraphics[width=11cm]{HUBBLE_PNMFSVD.eps} \\
}
\includegraphics[width=11cm]{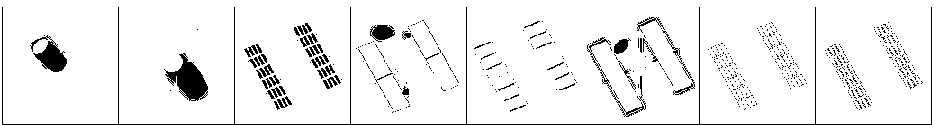} \\
\includegraphics[width=11cm]{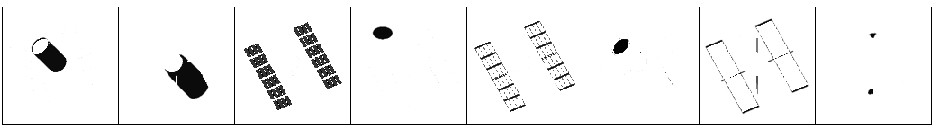}  \\
\end{tabular}
\caption{Hubble data set decomposition. From top to bottom: sample images at different wavelengths along with the true constituent materials, $k$-means, spherical $k$-means, CHNMF, CH(SVD), \ProjNMF, \ProjNMFSVD, EM-ONMF and ONP-MF.}
\label{HubbleResults_rnd}
\end{figure*}

We first use a synthetic data set from \cite{Hubble}, see Figure~\ref{HubbleResults_rnd} (top row), in clean conditions (i.e., without noise or blur). 
It represents the Hubble telescope and is made up of 8 different materials, each having a specific spectral signature. 
Figure~\ref{HubbleResults_rnd} displays the clustering obtained by the different algorithms\footnote{For EM-ONMF, $k$-means and SKM we preprocess the data by discarding pixels from the background (i.e., all columns of the input matrix with zero $\ell_2$-norm). Recall that, for each algorithm, we keep the best solution (w.r.t.\@ the error) among the 30 randomly generated initial matrices. 
} 
and we observe that only ONP-MF is able to successfully recover all eight materials without any mixing. 
Even with the SVD-based initialization, CHNMF and 0-PNMF (i.e., CH(SVD) and  \ProjNMFSVD) are not able to separate all materials properly; ONP-MF is the only algorithm able to perform this task almost perfectly.

\subsubsection{Urban Data Set}

\begin{figure*}[ht!] 
\centering
\begin{tabular}{c}
\includegraphics[width=12cm]{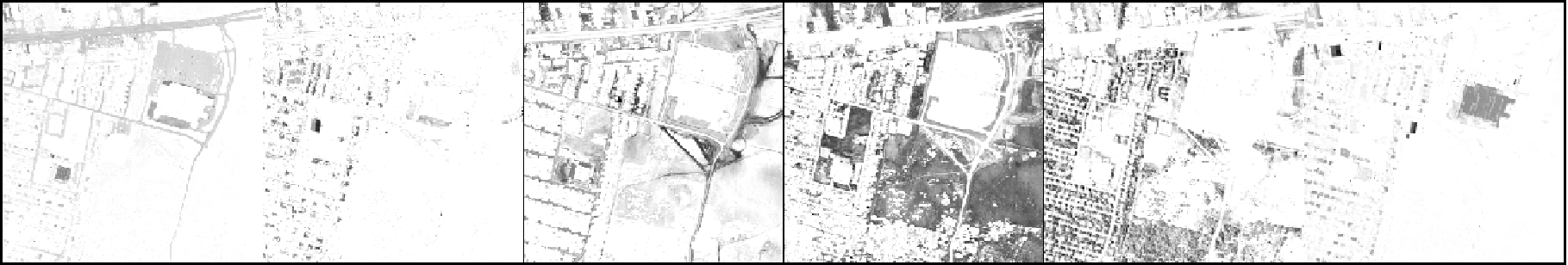} \\
\includegraphics[width=12cm]{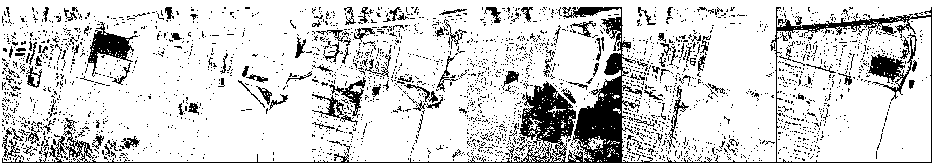} \\
\includegraphics[width=12cm]{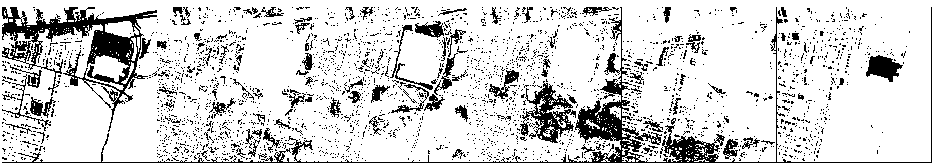} \\
\includegraphics[width=12cm]{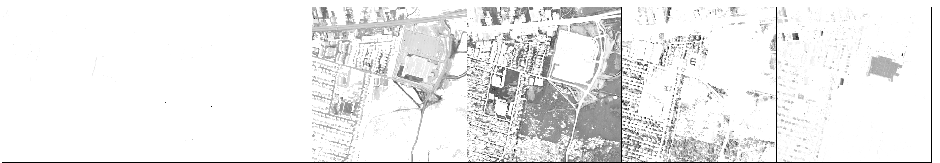} \\
\iftoggle{use_ORTH_PNMF_instead_of_PNMF}
{
\includegraphics[width=12cm]{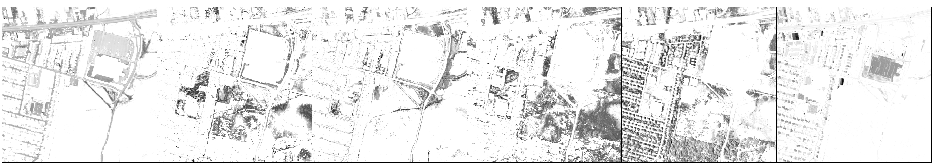} \\
}
{
\includegraphics[width=12cm]{URBAN_PNMFRAND.eps} \\
}
\includegraphics[width=12cm]{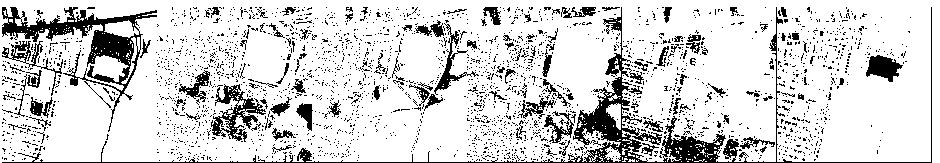} \\
\includegraphics[width=12cm]{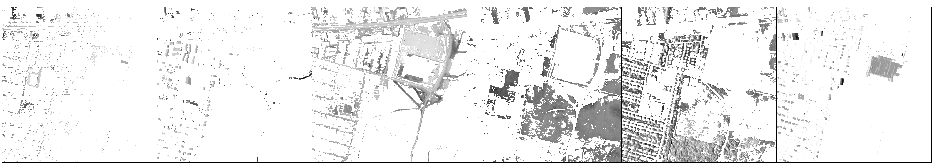}  \\
\end{tabular}
\caption{Urban data set decomposition. From top to bottom: `true' materials, $k$-means, spherical $k$-means, CHNMF, \ProjNMF, EM-ONMF and ONP-MF.}
\label{UrbanResults_rnd}
\end{figure*}

\begin{figure*}[ht!] 
\centering
\begin{tabular}{c}
\includegraphics[width=12cm]{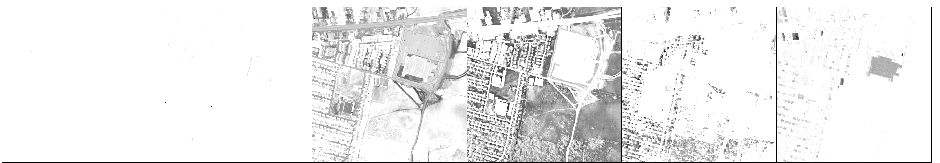} \\
\iftoggle{use_ORTH_PNMF_instead_of_PNMF}
{
\includegraphics[width=12cm]{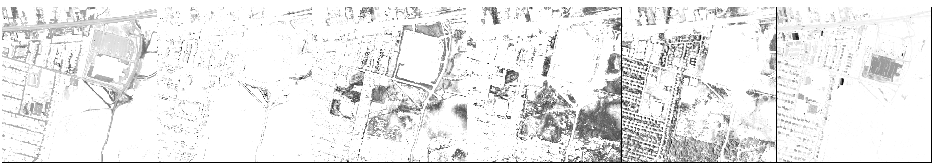} \\
}
{
\includegraphics[width=12cm]{URBAN_PNMFSVD.eps} \\
}
\end{tabular}
\caption{Urban data set decomposition. From top to bottom: CH(SVD) and \ProjNMFSVD.}
\label{UrbanResults_svd}
\end{figure*}

The Urban hyperspectral image is taken from  HYper-spectral Digital Imagery Collection Experiment (HYDICE) air-borne sensors. It contains 162 clean bands, and $307 \times 307$ pixels for each spectral image; it is mainly composed  of 6 types of materials: road, dirt, trees, roofs, grass and metal (mostly metallic rooftops) as reported in \cite{GWO09, GP10}. The first row of Figure~\ref{UrbanResults_rnd} displays a very good clustering obtained using N-FINDR5 \cite{Win99} plus manual adjustment from \cite{GWO09}, along with the clusterings obtained with the different algorithms. 
The road and dirt are difficult to extract because their spectral signatures are similar (up to a multiplicative factor), and none of the algorithms is able to separate them perfectly. 
ONP-MF successfully extracts the grass, trees, and roofs and is the only algorithm able to extract the metal (second basis element), while only mixing the road and dirt together. Spherical $k$-means, CHNMF, \ProjNMFSVD \, (Figure~\ref{UrbanResults_svd}) and EM-ONMF also perform relatively well, being able to extract the road (mixed with dirt or metal), trees, grass (as two separate basis elements) and roofs. 
CH(SVD) and $k$-means perform relatively poorly: they are not able to separate as many materials properly. 


\subsection{Image Segmentation: Swimmer Data Set} \label{ima}

\begin{figure*}[ht!]
\centering
\begin{tabular}{c}
\includegraphics[width=10cm]{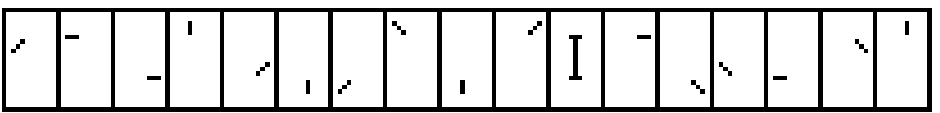} \\
\includegraphics[width=10cm]{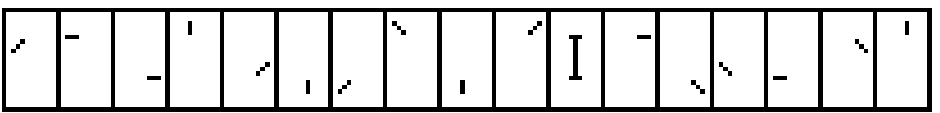} \\
\includegraphics[width=10cm]{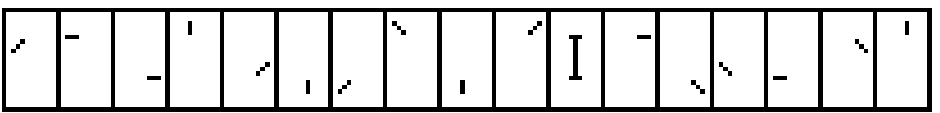} \\
\includegraphics[width=10cm]{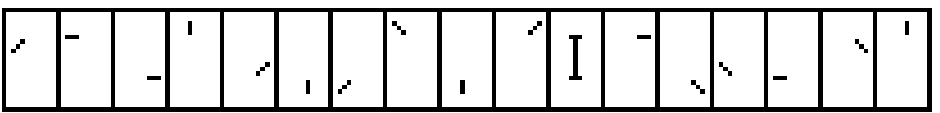} \\
\iftoggle{use_ORTH_PNMF_instead_of_PNMF}
{
\includegraphics[width=10cm]{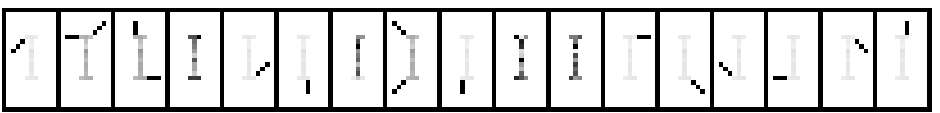}  \\
\includegraphics[width=10cm]{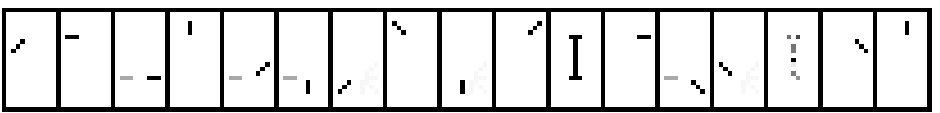}  \\
}
{
\includegraphics[width=10cm]{SWIMMER_PNMFRAND.eps}  \\
\includegraphics[width=10cm]{SWIMMER_PNMFSVD.eps}  \\
}
\includegraphics[width=10cm]{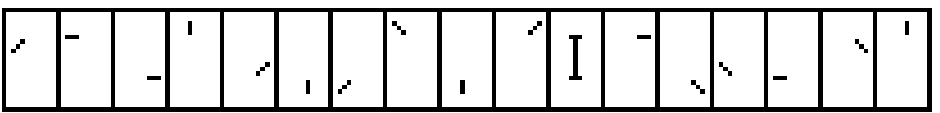}  \\
\includegraphics[width=10cm]{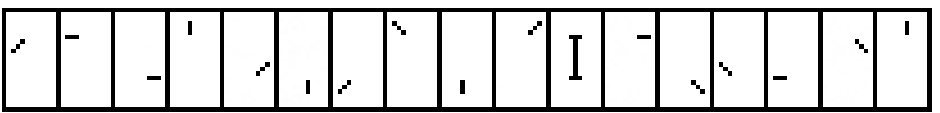}  \\
\end{tabular}
\caption{Swimmer data set decomposition. From top to bottom: $k$-means, spherical $k$-means, CHNMF, CH(SVD), \ProjNMF, \ProjNMFSVD, EM-ONMF and ONP-MF.}
\label{Swim_rnd}
\end{figure*}

The swimmer image data set consists of 256 binary images of a body with four limbs which can be positioned in four different ways each. The goal is to find a part-based decomposition of these images, i.e., isolate the different constitutive parts of the images (the body and the limbs, 17 in total). Moreover, these parts are not overlapping, and therefore no rows of $V$ can share nonzero entries in the same column, and ONMF is an appropriate model. 
Figure~\ref{Swim_rnd} displays  the basis elements obtained with the different ONMF algorithms. It can be observed that, in this case, the SVD-based initialization is of no benefit, neither for CHNMF nor for \ProjNMF. All algorithms are able to successfully find the correct parts except \ProjNMF \, and \ProjNMFSVD. 

\section{Conclusion} \label{con} 

In this paper, we have studied the ONMF problem and showed its equivalence with a weighted variant of spherical $k$-means (Theorem~\ref{Th1}). This led us to design a new EM-like algorithm for solving ONMF problems (Alg.~\ref{algomod}, EM-ONMF). We have also proposed an alternative approach based on an augmented Lagrangian method imposing orthogonality at each step while relaxing the nonnegativity constraint (Alg.~\ref{ONPMF}, ONP-MF). 

We then performed numerical experiments on some synthetic, text and image data sets. Note that Euclidean-based metric ONMF~\eqref{eq:JoNMF} is not particularly suited for document classification: First, the use of the Frobenius norm makes the implicit assumption that the noise is Gaussian, which is unrealistic for sparse data sets such as document data sets; see, e.g., the discussion in \cite{CK12}. Second, ONMF assumes that each document is on a single topic, which is in general not true (there exist more general generative models such as LDA assuming that documents are mixtures of several topics). Text data sets are nevertheless a worthy benchmark on which to compare the effectiveness of various Euclidean-metric ONMF algorithms.

The experiments indicate that our ONP-MF algorithm is by far the most robust among existing algorithms for solving the Euclidean ONMF problem~\eqref{eq:ONMF}: it always gave very good results, the best in many cases, using only one initialization. In particular, we observed for all image experiments that a single (deterministic) run of ONP-MF worked better than all the other tested algorithms, despite the fact that those were allowed to keep the best solution obtained from 30 different (random) initializations.  Since initialization is known to be an important component in the design of successful NMF methods~\cite{Bou08SVD}, we believe that initializing the $V$ factor with the unaltered right singular vectors of the data matrix, which is allowed by the workings of ONP-MF but impossible with other ONMF methods, plays an instrumental role in the clustering performance of ONP-MF observed in numerical experiments.

 \section*{Acknowledgments}
 
  The authors would like to thank the guest editor and the reviewers for their feedback which helped improve the paper.

\bibliographystyle{siam}
\bibliography{onpmf}

\end{document}